\documentclass[10pt,a4paper,twosided]{amsart}

\usepackage{amsmath}
\usepackage{amsfonts}

\theoremstyle{plain}
\newtheorem{theorem}{Theorem}[section]
\newtheorem{lemma}[theorem]{Lemma}
\newtheorem{proposition}[theorem]{Proposition}
\newtheorem{cor}[theorem]{Corollary}

\newtheorem{bigthm}{Theorem}

\newtheorem{mt}{Main Theorem}

\theoremstyle{definition}

\theoremstyle{remark}
\newtheorem{rem}[theorem]{Remark}

\numberwithin{equation}{section}

\newcommand{\zz}{{\mathbb{Z}}}
\newcommand{\cz}{{\mathbb{C}}}
\newcommand{\kz}{{\mathbb{K}}}
\newcommand{\pz}{{\mathbb{P}}}

\newcommand{\ra}{\rightarrow}
\newcommand{\cF}{\mathcal {F}}
\newcommand{\cE}{\mathcal {E}}
\newcommand{\cG}{\mathcal {G}}
\newcommand{\cL}{\mathcal {L}}
\newcommand{\cO}{\mathcal {O}}
\newcommand{\cI}{\mathcal {I}}
\newcommand{\rxa}[1]{\xrightarrow{\kern 1pt #1\kern 1pt}}
\newcommand{\lxa}[1]{\xleftarrow{\kern 1.5pt #1\kern 1pt}}

\DeclareMathOperator{\rk}{rk}
\DeclareMathOperator{\im}{im}
\DeclareMathOperator{\coker}{coker}
\DeclareMathOperator{\depth}{depth}

\begin{document}
\title{An explicit construction of ruled surfaces}
\author{Alberto Alzati}
\address{Dipartimento di Matematica Univ. di Milano\\
via C.\ Saldini 50 20133-Milano (Italy)}
\email{alzati@mat.unimi.it}
\thanks{This work is within the framework of the national research project ``Geomety
of Algebraic Varieties'' Cofin 2004 of MIUR 
and the DFG-Forschungsschwerpunkt ``Globale Methoden in der komplexen Geometrie''}
\author{Fabio Tonoli}
\address{Matematisches Institut der Universit\"at Bayreuth -- Lehrstuhl VIII\\
Universit\"atstra\ss e 30, D-95447 Bayreuth (Deutschland)}
\email{Fabio.Tonoli@uni-bayreuth.de}

\maketitle

\begin{abstract}
The main goal of  this paper  is  to give a general method
to  compute (via  computer algebra systems) 
an explicit set of generators of the ideals of the projective  
embeddings of some ruled surfaces, 
namely projective line bundles over curves such that  the  fibres are embedded
as  smooth rational curves.

Indeed, although the existence of the embeddings that we  consider
is well known, often in literature there are no explicit descriptions
of the corresponding projective ideals. 
Such an explicit description allows to compute,
besides all the syzygies, 
some of the important  algebraic invariants of the surface, 
for instance the $k$-regularity, which are not
always easy to compute by general formulae  or by geometric arguments.

An implementation of our algorithms and explicit examples for the computer
algebra system {\tt Macaulay2} (cf. \cite{GS}) are included,
so that anyone can use them for his own purposes.
\end{abstract}

\section*{Introduction and Notation}

Let $E$ be a rank $2$ vector bundle over a smooth, genus $q,$ curve $C$. 
It is known that any such vector bundle $E$, regarded as a sheaf, 
is an extension of invertible sheaves. 
If $E$ is a normalized vector bundle, i.e. $H^{0}(C,E)\neq 0$ 
but $H^{0}(C,E\otimes G)=0$ for any line bundle $G$ of negative degree,
then $E$ sits into a short exact sequence
\begin{equation}\label{E=extension}
0\rightarrow \mathcal{O}_{C}\rightarrow E \rightarrow L\rightarrow 0,
\end{equation}
and $L=\det (E)$.

We now consider the geometrically ruled surface $X:=\pz(E)$,
endowed with the natural projection $p:\ \pz(E)\ra C$.
In this case $Pic(X)\cong \zz\oplus p^{*}Pic(C)$, where $\zz$ is generated 
by the tautological divisor of $X$, i.e. a divisor $C_0$, image of a section 
$\sigma_0: C \ra X$, whose associated invertible sheaf is $\cO_X(1)$.
According to this notation, every divisor on $X$ is linearly
(resp. numerically) equivalent to 
$aC_{0}+p^{*}B$ (resp. $aC_{0}+bf$, being $f$ a fiber of $p$)
where $B$ is a degree $b$ divisor of $C$.

We choose a very ample divisor $A$ on $X$ and 
we consider the polarized ruled surface $(X,A)$:
what are then the equations of $X$?
In other words, $X$ is embedded in $\pz^{h^{0}(X,A)-1}$ by $|A|$ 
and we aim to give an algorithm for computing
a set of generators of the ideal $I_{X}$ in the
ring $S(V):=\oplus_{i\geq 0} S^i(V)$, the symmetric algebra of
$V=H^{0}(X,A)$.

Ampleness conditions for the divisor $A$ are classical and well known
(cf. e.g. \cite{h}).
In particular, by the Nakai's criterion, denoting with $e:=-\deg E$ the
invariant of $X$, an ample divisor $A$ is numerically equivalent to
$aC_{0}+bf$ with $a\geq 1$ and $b>ae$ if $e\geq 0$
or $b>\frac{ae}{2}$ if $e<0$, and the very ampleness for $A$ 
should be checked case by case with some criteria,
e.g. Reider's criterion (cf. \cite{rei}) or by looking at the image of
$X$ by $|A|$.

\medskip
The search of an algorithm has a positive answer. 
However the implementation of the algorithm obviously is not part of
the statement:
\begin{mt}
Let $C\subset \pz^m$ be a smooth curve $C$ of genus $q$,
$B$ a divisor on $C$ and $L$ a line bundle over $C$.
Consider a normalized rank 2 vector bundle $E\in\mathcal Ext^1(L,\cO_C)$ 
over $C$ given by an extension $0\ra \cO_C \ra E \ra L \ra 0$ and
suppose that the divisor $A=kC_{0}+p^{*}B $ on the surface
$X=\pz(E)$ is very ample. 
Then there is an algorithm yielding a set of generators for the ideal $I_{X}$ 
of the embedded $X$ in $\pz^{h^{0}(X,A)-1}=\pz(H^{0}(X,A)^{*})$ by $|A|$.
\end{mt}

\smallskip
The first step in the algorithm is to choose a projective model
for the curve $C$.  
This is not known in general, but an algorithm 
to pick curves up to genus 14 at random is implemented in \cite{ST},
and we refer to this algorithm for this part of the construction.

The second step is the choice of $E$. 
To this purpose we first give $L$ 
by giving explicitly a meromorphic section, i.e. writing
$L=\cO(D_1-D_2)$, where $D_1$ and $D_2$ are effective divisors on $C$,
and then we choose explicitly an extension, according to
Lemma \ref{mapcone}.

The final choice is the choice of $A$, determined by
the value of $k$ and another divisor $B$ on the curve $C$.
In order to perform the construction, another suitable divisor $D$ on $C$ has
to be chosen, but the polarized surface $(X,A)$ 
is independent from this choice.

\medskip
There are other ways to construct a polarized ruled surface $X$,
the most powerful being the one of considering
a locally free resolution of its ideal sheaf $\cI_X$.

Indeed, if $X$ has codimension 2, i.e. $X\subset \pz^4$, 
then there exist two sheaves $\cF$ and $\cG$  with
$\rk \cG=\rk F+1$ and a map
$\Phi:\ \cF \ra \cG$ such that the Eagon-Northcott complex
defined by the minors of $\Phi$ identifies $\coker \Phi$ with
a suitably twisted ideal sheaf of $X$.
The sheaves $\cF$ and $\cG$ are then constructed starting from
the cohomology table of $\cI_X$. This constructing method
was introduced in \cite{des} and largely used to
construct surfaces in $\pz^4$ (c.f. also \cite{ds} for 
a further description and a nearly up-to-date list of references).
If instead $X$ has codimension 3, this type of construction
can be still performed using the Pfaffian complex instead
of the Eagon-Northcott complex:
indeed $X$ is a codimension 3 subcanonical scheme in $\pz^5$
and a locally free resolution of its ideal sheaf is still
known (c.f.  \cite{w}).

We hope that our method will enable 
to construct new examples of polarized ruled surfaces whose existence
is not known, because it allows more control 
over the geometry of the polarized surface (we mean: the choices of
the curve $C$, of the extension giving $E$, and of the divisor $A$).

\medskip
We remark also that such an algorithm for a particular case of scrolls
(ruled surfaces over a curve of genus $2$ embedded
as scrolls of degree 8 in $\pz^5$) 
was partially described in \cite{c}.

\medskip
The paper is structured as follows.
The first section is devoted to the construction of scrolls, which
is the core of the construction. In particular, we choose a divisor
$D$ in order to construct the module $H^0_*(C,E\otimes\cO_C(B))$:
once this is done the algorithm is straightforward.
The third section treats the case of conic bundles, while the
fifth section generalizes the algorithm for conic bundles to
the case of $k$-bundles.

The even sections instead contain the implementation
of the algorithm with the computer--algebra program {\tt Macaulay2}
(\cite{GS}) and examples for each case:
a family of scrolls of degree 8 in $\pz^5$ with sectional genus 2
(cf. \cite{c}),
a family of scrolls of degree 6 in $\pz^5$ with sectional genus 1,
and two families of ruled surfaces in $\pz^5$: 
conic bundles of degree 8 with sectional genus 3
and cubic bundles of degree 9 with sectional genus 4.

\medskip
\begin{center} \bf Notation Table:\end{center}
{
\small 
\newcommand{\linestretch}{\vspace{-2pt}}
\begin{tabular}{lp{10cm}}
\linestretch
$\kz$ &base field, usually $\cz$\\
\linestretch
$\pz^{n}$ &projective $n$-dimensional space over $\kz$\\
\linestretch
$\pz(E)$ &projectivization of the rank $2$ vector bundle $E$ over a
smooth curve $C$, 
$C_{0}$ is its tautological divisor,
$p:\pz(E)\rightarrow C$ the natural projection, and
$f$ the numerical class of a generic fibre of $p$\\
\linestretch
$c_{i}(E)$ &i-th Chern class of $E$\\
\linestretch
$\mathbf{F}_{e,q}$ &ruled surface of invariant $e:=-\deg
[c_{1}(E)]\geq -q$ over a smooth, genus $q,$ curve $C$\\
\linestretch
$\equiv $ $$ &numerical equivalence\\
\linestretch
$^{*}$ &means duality\\
\linestretch
$|D|$ &linear system of effective divisors linearly equivalent to the
divisor $D$\\
\linestretch
$I_W$ ($\cI_W$)  & ideal (ideal sheaf) of a projective variety 
$W\subset \pz^{n}$\\
\linestretch
$K_{W}$ &canonical divisor of a smooth variety $W$\\
\linestretch
$g(W)$ &sectional genus of a smooth variety $W\subset \pz^{n}$\\
\linestretch
$H_{*}^{0}(W,\mathcal{F})$ &$\bigoplus_{t\geq 0}H^{0}(W,\mathcal{F}\otimes 
\mathcal{O}_{W}(t))$ for any sheaf $\mathcal{F}$ on $W\subset \pz^{n}$\\
\linestretch
$\widetilde{M}$ &sheaf of $\mathcal{O}_{W}$-modules associated to any 
$S$-module $M$, where $S$ is the coordinate ring of a smooth variety $W$\\
\linestretch
$S(V)$  &$\oplus_{n\geq 0}S^n(V)$ symmetric algebra of the vector space $V$\\
\linestretch
$S(E)$ &$\oplus_{n\geq 0}S^n(E)$ symmetric $\cO_W$-algebra 
of the vector bundle $E$ over a variety $W$\\
\linestretch
$\mu (E)$ &$\deg E/\rk E $, slope of the vector bundle  $E$\\
\linestretch
$\mu ^{-}(E)$ &$\min \{\mu (Q)|E\rightarrow Q\rightarrow 0\}$
\end{tabular}
}

\section{Construction of scrolls}

In this section we give an algorithm to compute explicity a set of generators
for the ideal of embedded scroll surfaces.

\begin{bigthm}
\label{scrolls}Let $C\subset \pz^m$ be a smooth curve $C$ of genus $q$,
$B$ a divisor on $C$ and $L$ a line bundle over $C$.
Consider a normalized rank $2$ vector bundle $E\in\mathcal Ext^1(L,\cO_C)$ 
over $C$ given as extension $0\ra \cO_C \ra E \ra L\ra 0$
and suppose that the divisor $A=C_{0}+p^{*}B$ on the surface
$X=\pz(E)$ is very ample.
Then there is an algorithm yielding 
a set of generators for the ideal $I_{X}$ 
of the embedded $X$ in $\pz^{h^{0}(X,A)-1}=\pz(H^0(X,A)^*)$ by $|A|$.
\end{bigthm}

Before giving the proof of the theorem, in terms of an explicit algorithm, let
us first point out some remarks and technical Lemmas.
Let $I_{C}$ be the ideal of the curve $C$ in $\pz^{m}=Proj(R)$, 
where $R=\kz[x_{0},x_{1},...,x_{m}]$, and let $S:=R/I_{C}$ be the 
coordinate ring of $C\subset \pz^{m}$.

At first, we remark here that the algorithm is a straigthforward computation once
we present the module  
$M$ defined as 
$$M:=H^0_*(C,E \otimes\cO_{C}(B))=\oplus_{i\geq 0}H^0(C,E \otimes\cO_{C}(B+iH)),$$
where $H$ is an hyperplane divisor of $\pz^{m}$.
The details of this computation will be given later in the proof of the theorem.
The short exact sequence  (\ref{E=extension}) implies the exactness of
\begin{equation}\label{EB}
0\ra \cO_{C}(B)\ra E \otimes\cO_{C}(B) \ra L\otimes\cO_{C}(B)\ra 0,
\end{equation}
from which we can derive the desired presentation of $M$.

\medskip
Next, we present here some Lemmas needed for the proof of the theorem.

\begin{lemma}
\label{moduli}Let $I_{C}$ be the ideal of a smooth curve $C$ in a projective
space $\pz^{m}$ and let $S$ be the coordinate ring of $C$.
Let $D$ be an effective divisor on $C$. Then
the $S$-modules $H_{*}^{0}(C,\mathcal{O}_{C}(D))$ and $
(I_{D})^{*}:=Hom_{S}(I_{D},S)$ are naturally isomorphic as (graded) $S$-modules, 
where $I_{D}\subset S$ is the ideal of the divisor $D.$
\end{lemma}

\begin{pf}
Let us recall the following well-known result on local cohomology 
(cf. \cite[Thm. A4.1]{e}).
Let $S$ be a graded noetherian ring with degree $0$ part a field, 
$\mathfrak m$ a maximal ideal and $M$ a finitely generated $S$-module. 
Then there is a natural exact sequence
$$0 \ra H_{\mathfrak m}^0(M) \ra M \ra \oplus_{i\geq 0} H^0({\rm
Proj}\ S,\widetilde M(i)) \ra  H_{\mathfrak m}^1(M) \ra 0,$$
where $H_{\mathfrak m}^i(M)$ denotes the $i$-th local cohomology group of $M$
with respect to $\mathfrak m$.
If you consider a variety $W\subset $ $\pz^{m}$ and 
its coordinate ring $S$ and you take any graded $S$-module $M,$ 
then $M=H_{*}^{0}(W,\widetilde{M})$ 
if $H_{\mathfrak{m}}^{0}(M)=H_{\mathfrak{m}}^{1}(M)=0,$ where $\mathfrak{m}$ is the
maximal ideal of $S$.
Hence $M\cong H^0_*(\widetilde M)$
if  $H_{\mathfrak{m}}^{0}(M)=H_{\mathfrak{m}}^{1}(M)=0$.
Moreover the two vanishings $H_{\mathfrak{m}}^{0}(M)=H_{\mathfrak{m}}^{1}(M)=0$ 
follow from the condition $depth(\mathfrak{m},M)\geq 2$ by \cite[Ex. 3.4 and 3.3]{h}.
For more details and related results, see \cite{g}
(in particular Prop. 2.2 and Thm. 3.8).

In our case $W=C$, 
$S:=\kz[x_0,\dots,x_m]/I_C$ is the coordinate ring of the curve $C$ in $\pz^{m}$, 
$\mathfrak m:=(x_0,\dots,x_m)$ is the image of the irrelevant ideal of $\kz[x_0,\dots,x_m]$
in $S$, and $M=(I_{D})^{*}$. 
Remark that $S$ is the coordinate ring of a cone over a curve,
and therefore $I_D$, as well as $(I_D)^*$, is not necessarily a projective
$S$-module, since it may not be locally free
in the local ring at the vertex of the cone.
Therefore we proceed as follows.

At first, notice that $\depth(\mathfrak m,S)=2$ and that a regular sequence
for $S$ is also a regular sequence for $I_D$,
since $I_D$ is a submodule of $S$.
In second place, remark the following:
if $t_1,\dots,t_d\in\mathfrak m$ is a regular sequence for $S$, then it is also a 
regular sequence for $(I_{D})^{*}=Hom(I_D,S)$.
We argue by contradiction.
Suppose that $t_i$ is a $0$-divisor for 
$(I_{D})^{*} \mod{(t_1,\dots,t_{i-1})}$. 
Then there esists a non-zero morphism
$\varphi\in (I_{D})^{*} \mod{(t_1,\dots,t_{i-1})}$ s.t. 
$t_i\varphi=0 \mod{(t_1,\dots,t_{i-1})}$.
Take an $x\in I_D$ s.t. $\varphi(x)\neq0$ in $S/{(t_1,\dots,t_{i-1})}$:  
from $t_i\varphi(x)=0$ in  $S/{(t_1,\dots,t_{i-1})}$ we get 
that $t_i$ is a $0$-divisor in $S/{(t_1,\dots,t_{i-1})}$, a contradiction.
Moreover $t_1,\dots,t_d\in\mathfrak m$, 
hence we have $(t_1,\dots,t_d)(I_{D})^{*}\neq (I_{D})^{*}$.
Indeed, if this is not the case, then $\mathfrak m
(I_{D})^{*}=(I_{D})^{*}$ and therefore there exists an element 
$r\in\mathfrak m$ such that $(1-r)(I_{D})^{*}=0$,
cf. \cite[Cor. 4.7]{e}.
In particular, considering the inclusion $\iota:\ I_D\ra S$,
we have $(1-r)\iota=0$ and therefore $1-r$ is a 0-divisor in $S$.
Since $S$ is an integral domain, it follows that $r=1$, which is absurd 
since $r\in\mathfrak m$.

We conclude that  
$\depth(\mathfrak m,(I_D)^*)\geq\depth(\mathfrak m,S)=2$
and therefore we get 
$(I_{D})^{*}=H_{*}^{0}(C,\widetilde{[(I_{D})^{*}]})=H_*^0(C,\cO_{C}(D))$, 
as well as 
$\depth(\mathfrak m,I_D)\geq\depth(\mathfrak m,S)=2$ and $I_{D}=
H_*^0(C,\cO_{C}(-D))$.
\end{pf}

\begin{lemma}
\label{mapcone} Let $F$ and $G$ be two $S$-modules with free
resolutions:
$$F^\bullet : \qquad \cdots \ra F_3 \rxa{\phi_3} F_2 
\rxa{\phi_2} F_1\rxa{\phi_1} F_0 \rxa{\phi} F\ra 0,$$
$$G^\bullet : \qquad \cdots \ra G_3 \rxa{\psi_3} G_2 
\rxa{\psi_2} G_1\rxa{\psi_1} G_0 \rxa{\psi} G\ra 0.$$
Then any morphism $\varphi\in Hom_S(F_1,G_0)$ satisfying
$\psi\circ\varphi\circ\phi_2=0$, i.e. inducing a morphism in
$Hom_S(\ker\phi,G)=Hom_S(\im\phi_1,\coker\psi_1)$, determines an extension $M\in Ext^1_S(F,G)$
and, conversely, any extension is determined by such a morphism.

Moreover, the module $M\in Ext^1_S(F,G)$ corresponding 
to $\varphi$ has presentation 
$$\begin{pmatrix}\phi _{1} & 0 \\ \varphi  & \psi _{1}\end{pmatrix} 
:\ F_1\oplus G_1 \ra F_0\oplus G_0.$$
\end{lemma}

\begin{pf}
Consider the module $K=\ker\phi$.
The short exact sequence $0 \ra K \ra F_0 \rxa{\phi} F \ra 0$
induces by duality
$$Hom_S(F_0,G) \ra Hom_S(K,G) \ra Ext_S^1(F,G) \ra Ext_S^1(F_0,G)=0.$$
Therefore $ Ext_S^1(F,G)\cong Hom_S(K,G)/Hom_S(F_0,G)$.

In the same way the short exact sequence 
$0\ra \ker\phi_1 \ra F_1 \rxa{\phi_1} K\ra 0$ gives
$0\ra Hom_S(K,G)\ra Hom_S(F_1,G)\ra Hom_S(\ker\phi_1,G)$.
Hence $Hom_S(K,G)$ is the kernel of the second map and we can identify $Hom_S(K,G)$ 
with the set of morphisms $\eta\in Hom_S(F_1,G)$
whose restriction to $\ker\phi_1$ is zero or equivalently, 
since $\ker\phi_1=\im\phi_2$, whose composition $\eta\circ\phi_2$ is zero. 
Since $F_1$ is projective, the surjectivity of $G_0 \ra G$ gives
the surjectivity of $Hom_S(F_1,G_0) \ra Hom_S(F_1,G)$.

In conclusion, a morphism $\varphi\in Hom_S(F_1,G_0)$ satisfying the
hypothesis  determines by composition a morphism in
$\eta\in Hom_S(F_1,G)$ satisfying $\eta\circ\phi_2=0$
which therefore is an element of $Hom_S(K,G)$:
its rest class in $Hom_S(K,G)/Hom_S(F_0,G)$ determines an extension 
$M\in Ext_S^1(F,G)$, as desired.

To compute a presentation of such an extension $M$, let us
denote with $\iota$ the inclusion $K\ra F_0$, and with $\varphi'\in Hom_S(K,G)$ 
the  morphism induced by $\varphi$.
Then the module $M$ is the quotient $(F_0\oplus G)/\im (\iota\oplus \varphi')$,
cf. \cite[pag. 722]{gh} or \cite[Ex. A3.26]{e},
which is the cokernel of the morphism 
$\begin{pmatrix}\phi _{1} & 0 \\ \varphi  & \psi _{1}\end{pmatrix} 
:\ F_1\oplus G_1 \ra F_0\oplus G_0.$
\end{pf}

\begin{lemma}
\label{risoluzione} 
Let $W\subset \pz^m$ be a smooth algebraic variety and let
$\cE$ be a locally free sheaf on $W$. 
Denote by $S:=\kz[x_0,\dots,x_m]/I_W$ the coordinate ring of $W$ in $\pz^m$.
Suppose further that the tautological bundle
$\tau_{\pz(\cE)}=\cO_{\pz(\cE)}(1)$ of $\pz(\cE)$ is very ample.
Then, given a presentation of the $S$-module $M:=H^0_*(W,\cE)$,
there is an algorithm yielding a set of generators for the ideal 
$I_{\pz(\cE)}$ of the embedded variety $\pz(\cE)$ 
by the complete linear system $H^0(\pz(\cE),\tau_{\pz(\cE)})=H^0(W,\cE)$.
\end{lemma}

\begin{pf}
Let $h^0(W,\cE)=n+1$ and let $\kz$ be the base field of $W$.
The given embedding $\iota$ associated 
to the complete linear system $H^0(\pz(\cE),\tau_{\pz(\cE)})=H^0(W,\cE)$
comes with a map of sheaves of rings on $\pz^n$
$\iota^\#:\ \cO_{\pz^n} \ra i_*\cO_{\pz(\cE)}(1)$,
induced by sending $n+1$ new variables $y_0,\dots,y_n$ to the global
sections of $H^0(W,\cE)$, which generate the $\cO_W$-algebra
$S(\cE)=\oplus_{d\geq0} S^d(\cE)$.
The ideal sheaf $\widetilde I_{\pz(\cE)}$ is given by the kernel of this map.

Let $M'\subset M$ be the $S$-submodule generated by a basis of
$H^0(W,\cE)$.
If $\phi$ is the given free presentation of $M$,
we can compute a free presentation $\phi'$ of $M'$ of the form:
$$ M_1 \rxa{\phi'} M_0 \ra M'\ra 0,$$
where $\rk M_0=n+1$, i.e. the generators of $M_0$ map to a base of
$H^0(W,\cE)$
(the required algorithm computes the relations among the given set of
generators of $M'$, and it is usually implemented in computer-algebra
programs).

Consider in $S[y_0,\dots,y_n]$ the ideal $I$ given by 
\begin{equation} \label{RelProj}
I:=\begin{pmatrix}  y_0 &\dots &y_n\end{pmatrix} \cdot \phi'.
\end{equation}

The ideal $I_{\pz(\cE)}$ is given by the polynomial relations
among the $\{y_{0},\dots,y_{n}\}$ in the saturation of $I$ 
with respect to the ideal $(x_{0},...,x_{m})\subset S[y_0,\dots,y_n]$. 
Therefore $I_{\pz(\cE)}$ can be obtained by saturating $I$
with respect to the ideal $(x_{0},...,x_{m})$
and intersecting this new ideal with the subring $K[y_0,\dots,y_n]$.
\end{pf}

\begin{rem}
Equation (\ref{RelProj}) yields  a presentation, as $S[y_0,\dots,y_n]$-module, 
of the $S$-algebra generated by $H^0(W,\cE)$ in $S(H^0_*(W,\cE))$.
If $M$ is generated by $H^0(W,\cE)$, then $M$ admits a presentation 
$\oplus_{j=0}^s S(-l_j) \rxa{\phi} \oplus_{i=0}^n S \ra M \ra 0$
and the $S$-algebra $S(H^0_*(W,\cE))$ has a presentation 
$$\oplus_{j=0}^s S[y_0,\dots,y_n](-l_j) 
\rxa{(\dots,\sum_{i=0}^n y_i\phi_{ij},\dots)}
S[y_0,\dots,y_n] \ra S(H^0_*(W,\cE))\ra 0.$$
If this is not the case, i.e. $M$ is not generated by the minimal
degree part, let 
$\oplus_{j=0}^s S(-l_j)$ $ \rxa{\phi} \oplus_{i=0}^n S(-h_i) \ra M \ra 0$
be a presentation of $M$.
Then the $S$-algebra $S(H^0_*(W,\cE))$ has still a presentation as above,
but now the the ring $S[y_0,\dots,y_n]$ is weighted,
$y_i$ having weight $h_i$.
\end{rem}

\begin{rem}\label{m=m'}
Let us assume that $W = C\subset \pz^m$ is a smooth projective curve 
of genus $q$ and $\cE$ is a vector bundle over $C$.
If $\mu^-(\cE)\geq 2q$ and $\deg(\cO_C(1))\geq 2q$ and moreover one of
these two inequalities is strict, 
then $H^0(C,\cE)$ generates $H^0_*(C,\cE)$ as $S$-module,
where $S$ denotes the coordinate ring of $C$ in $\pz^m$.
\end{rem}
\begin{pf}
This is a direct application of Theorem 2.1 of \cite{bu}:
under these hypothesis the map
$H^0(C,E)\otimes H^0(C,\cO_C(t)) \ra H^0(C, E\otimes \cO_C(t))$ is
surjective $\forall t\geq 0$.
\end{pf}

\medskip
\begin{pf}[Proof of Thm. \ref{scrolls}]
Recall the assigned exact sequence in the statement
 $0\ra \cO_C \ra E \ra L\ra 0$ and 
let $D$ be an effective divisor on $C$ such that $D-B$ is effective satisfying 
both conditions
\begin{equation}
\label{direct constr scrolls}
\begin{cases}
H^1(C,\cO_C(D+jH))=0\text{ for }j\geq 0\\
|L\otimes \cO_C(D)|\text{ is not empty}
\end{cases},
\end{equation}
where $H$ is the divisor induced on $C$ by a hyperplane section of $\pz^m$. 

From the first condition of (\ref{direct constr scrolls}) 
there is a short exact sequence of cohomology modules 
$$0\ra H^0_*(C,\cO_C(D)) \ra H^0_*(C,E\otimes \cO_C(D)) \ra 
H^0_*(C,L\otimes \cO_C(D)) \ra 0,$$
implying that $H^0_*(C,E\otimes \cO_C(D))$ can be
obtained as an extension in 
$Ext_S^1(H^0_*(C,L\otimes \cO_C(D)),H^0_*(C,\cO_C(D)))$.\

The second condition of (\ref{direct constr scrolls}) implies 
the existence of an effective divisor
$D_2\in  H^0(C,L\otimes\cO_C(D))$.
Applying Lemma \ref{moduli} to the divisors $D_1=D$ and $D_2$, 
we get an explicit description of the modules $H^0_*(C,\cO_C(D))=I_{D_1}^*$ and
$H^0_*(C,L\otimes \cO_C(D))=I_{D_2}^*$.
In this way it is easy to compute their presentation in a computer algebra
system (cf. next section).

Given their presentations, since $H^0_*(C,E\otimes \cO_C(D))$ is
an extension in
$Ext_S^1(H^0_*(C,L\otimes \cO_C(D)), H^0_*(C,\cO_C(D)))$,
we apply Lemma \ref{mapcone} 
to get a presentation of $H^0_*(C,E\otimes \cO_C(D))$.
The tensorization of $H^0_*(C,E\otimes \cO_C(D))$ with
$I_{D-B}=H^0_*(C,\cO_C(-D+B))$ yields a module $M'$ 
whose associated coherent sheaf is $E\otimes B$.
The desired module $M=H^0_*(C,E \otimes\cO_{C}(B))$ can 
be obtained as $M=H^0_*(C,\widetilde M')$.
  
Finally, Lemma \ref{risoluzione} provides the requested set of
generators for $I_{X}$, where $X=\pz(E)$ is embedded by the 
very ample divisor $A=C_{0}+p^{*}B$.
\end{pf}

\begin{rem}\label{originale}
Suppose that $|B|$ and $|L\otimes B|$ contain effective divisors
$D_1\in |B|$ and $D_2\in|L\otimes B|$,
and that $\deg B>2q-2$ or $h^1(C,\cO_C(B+jH))=0\ \forall j\geq 0$.
Then in the proof of the theorem we can choose $D=B$, i.e.
the module $M=H^0_*(C,E \otimes\cO_{C}(B))$ can be directly obtained 
as extension in $Ext_S^{1}((I_{D_2})^*,(I_{D_1})^*)$.
\end{rem}

\begin{rem}
In the algorithm of the theorem, in order to get the whole module 
$M=H^0_*(C,E \otimes\cO_{C}(B))$ instead of just the submodule $M'$,
we use the corresponding implemented command in the computer algebra 
{\tt Macaulay2}.
If this command is not available in other computer systems, the algorithm in
the theorem is still valid under the assumptions in the Remark 
(\ref{originale}), assumptions which allow to choose $D=B$.
\end{rem}

\begin{rem}
The theorem is mainly used to obtain examples by considering 
random effective divisors $D_1,D_2$ with fixed degrees such that 
$\deg D_1>2q-2$
and a random extension class in $Ext_S^1((I_{D_2})^*,(I_{D_1})^*)$.
Defining $L$ as the sheaf $\cO_C(D_2-D_1)$, the previous extension
class determines one in $Ext_S^{1}(L,\cO_C)\cong H^1(C,L^*)$ and
the condition $\deg D_1>2q-2=0$ ensures that every extension 
in $Ext_S^{1}(L,\cO_C)$ can be obtained starting from an extension 
in $Ext_S^1((I_{D_2})^*,(I_{D_1})^*)$.
\end{rem}

\section{Some examples of ``interesting'' ruled surfaces}

In this section we will construct some examples 
of ``interesting'' ruled surfaces by applying Theorem \ref{scrolls},
where by ``interesting'' we mean that 
these surfaces have some particular properties.

\subsection{First example}
Let $C$ be a smooth curve of genus $2.$ Let $E$
be a normalized rank $2$ vector bundle of degree $2,$ so that we have the
following exact sequence:

\begin{center}
$0\rightarrow \mathcal{O}_{C}\rightarrow E$ $\rightarrow L\rightarrow 0$
\end{center}

where $L=\det (E)=c_{1}(E),$ $\deg L=2.$ Let $B$ be any degree $3$ divisor
of $C.$ On the surface $X=\pz(E)$ we can consider the divisor $%
A=C_{0}+p^{*}B\equiv C_{0}+3f.$ $A$ is a very ample divisor, whatever $B$ is
chosen, $h^{0}(X,A)=h^{0}(C,E\otimes B)=6$ and it embeds $X$ in $\pz^{5}$
as a smooth scroll of degree $8$ (see \cite{i2}); $g(X)=2.$ It is easy to
see that $X$ is $2$-normal if and only if it is not contained in a quadric.

About this surface we have the following proposition (see \cite{abb1}):

\begin{proposition}
\label{prop1}Let $X$ be the surface above. $X$ is contained in a rank $4$
quadric cone whose vertex is a $4$-secant line for $X$ and therefore $X$ is
not $2$-normal.
\end{proposition}

Note that in \cite{abb1} the proposition is proved by using geometric
arguments and it is not considered the $k$-normality of $X$ for $k\geq 3,$
moreover there is not a free resolution for the ideal $I_{X}.$ Some more
informations about $X$ can be found in \cite{c}, the article which suggested
to us to approach the problem.

\medskip
The choice of a random smooth genus 2 curve $C$  is performed by the following scripts
introduced in \cite{ST}:
\begin{small}
\begin{verbatim}
randomGenus2Curve = (R) -> (
     correctCodimAndDegree:=false;
     while not correctCodimAndDegree do (
          alpha:=transpose (vars R++vars R)**R^{-2} || 0*random(R^1,R^{2:-2});
          rd:=random(R^{8:-1,1:0},R^{6:-1});
          mappingCone:=rd|alpha;
          I:=ideal mingens ideal syz transpose mappingCone;
          correctCodimAndDegree=(codim I==2 and degree I==5););
     I);
isSmoothSpaceCurve = (I) -> (
     --I generates the ideal sheaf of a pure codim 2 scheme in P3
     singI:=I+minors(2,jacobian I);
     codim singI==4);
\end{verbatim}
\end{small}
The curve $C$ will be a degree 5 curve in $\pz^3$ (again, cf. \cite{ST}).
We also need to pick up  $t$ random points on $C$, which we perform by
separating the points of a good hyperplane section
(on a non algebraically closed field it can happen that these 
points are not separated):
\begin{small}
\begin{verbatim}
randomPoint = (C) -> (
     R:=ring C;
     isSinglePoint:=false;
     while not isSinglePoint do (
          hypsection:=C+ideal random(R^1,R^{-1});
          pt:=(decompose hypsection)#0;
          isSinglePoint=(degree pt==1););
     pt);
randomPoints = (C,t) -> (
     pt:=randomPoint C;i:=t-1;
     while i!=0 do (pti=randomPoint C;pt=intersect(pt,pti);i=i-1;);
     pt);
\end{verbatim}
\end{small}

\medskip
We give now the script. 
We choose the smooth genus 2 curve in $\pz^3$ with ideal
{\tt C}:
\begin{small}
\begin{verbatim}
K=ZZ/101;
R=K[x_0..x_3]
C=randomGenus2Curve R 
isSmoothSpaceCurve(C)
betti res C --it shows that the genus is really 2
\end{verbatim}
\end{small}
In this example we therefore set $m=3$ (the value $n=5$ is fixed), 
and we choose two effective divisors $L$ and $B$ of degree respectively 2 and 3:
in this way, by applying Remark \ref{originale}, we can settle $D=B$.
The points of $L$ and $D$ are chosen via the function  
{\tt randomPoints()},
which returns their ideals {\tt Ldual} and {\tt Ddual} in $\pz^3$.
Finally we compute the ideal in $\pz^3$ of the points of $L+D$, called {\tt D2dual}. 
\begin{small}
\begin{verbatim}
Ldual=randomPoints(C,2)
Ddual=randomPoints(C,3)
D2dual=intersect(Ldual,Ddual)
\end{verbatim}
\end{small}

We now compute the modules $H^0_*(L+D)$ and $H^0_*(D)$, called resp.
{\tt D2S} and {\tt DS}, 
where $S$ is the coordinate ring {\tt R/C} of $C$:
\begin{small}
\begin{verbatim}
S=R/C
DSdual=substitute(Ddual,S);DS=Hom(DSdual,S);
D2Sdual=substitute(D2dual,S);D2S=Hom(D2Sdual,S);
\end{verbatim}
\end{small}

We proceed, as explained in the mapping cone Lemma \ref{mapcone},
to compute a presentation {\tt phi} of a random module $M$ in
$Ext^1(H^0_*(L+D),H^0_*(D))$. 
We define for this purpose the function {\tt randomExt()} (cf. Lemma \ref{mapcone}):
\begin{small}
\begin{verbatim}
randomExt = (A,B) -> (
     phia:=presentation A;
     phib:=presentation B;
     Homom:=Hom(image phia,coker phib);
     phiab:=homomorphism random(Homom,S^1);phiab=matrix phiab;
     phiNull:=0*random(target phia,source phib);
     phi:=(phia||phiab)|(phiNull||phib);
     coker phi)
\end{verbatim}
\end{small}
and we apply this function to {\tt D2S} and {\tt DS}:
\begin{small}
\begin{verbatim}
M=randomExt(D2S,DS)
apply(-3..10,i->hilbertFunction(i,M))
\end{verbatim}
\end{small}
The module $M$ will be then the choice of 
$H^0_*(C,E\otimes \cO_C(D))$ corresponding to the choices of $L$, $D$
and the extension class in $Ext^1(H^0_*(L+D),H^0_*(D))$,
as explained in the proof of Theorem \ref{scrolls}.
The last line is a further (not needed) check on the Hilbert function of $M$.
We also remark that here $M$ is really $H^0_*(\widetilde M)$, 
because of our choice $D=B$.
Otherwise, the following line would compute the whole $H^0_*(\widetilde M)$:
\begin{small}
\begin{verbatim}
M=HH^0((sheaf M)(>=0));
\end{verbatim}
\end{small}

\medskip
We are now ready to compute an explicit set of generators of the ideal
$I_X$ of $X\subset\pz^5$, as explained in the proof of Lemma
\ref{risoluzione}. 
We define for this purpose the function {\tt scrollIdeal()}:
\begin{small}
\begin{verbatim}
scrollIdeal = (M) -> (
     phi=presentation prune image basis(0,M);
     T=K[y_1..y_(numgens target phi)];
     R:=ring phi;TR:=T**R;
     Phi:=substitute(phi,TR);
     IS:=ideal(substitute(vars T,TR)*Phi);
     J:=saturate(IS, ideal substitute(vars R,TR));
     ideal mingens substitute(J,T))
\end{verbatim}
\end{small}
As required in the proof of Lemma \ref{risoluzione}, 
the first line computes a presentation of the submodule generated 
by 
the elements of degree 0 of $M$,
i.e. by $H^0(C,E\otimes \cO_C(D))$, 
even if in this example this step is not needed by Remark \ref{m=m'},
since $\deg C=5$.

\medskip
We now  perform explicitely all the desired computations on 
$X\subset \pz^5$.
We firstly call this function and obtain $I_X$, called {\tt J} in the
script. Then  we check that $X$ is a smooth surface of degree 6, 
and we give the Betti table of $\cI_X$:
\begin{small}
\begin{verbatim}
J=scrollIdeal(M)
dim J, degree J
o26 = (3, 8)
codim (J+minors(3,jacobian J))
o27 = 6
betti res J
o28 = total: 1 8 15 13 6 1
          0: 1 .  .  . . .
          1: . 1  .  . . .
          2: . 6  7  . . .
          3: . 1  8 13 6 1
\end{verbatim}
\end{small}

By looking at the degrees of the set of generators for {\tt J},
it is easy to see that $X$ is contained in the following quadric cone
{\tt Q}, which is a rank $4$ quadric cone having the 4-secant line {\tt L}
as vertex, according to Proposition \ref{prop1}:
\begin{small}
\begin{verbatim}
Q=(gens J)_{0}
rank jacobian transpose jacobian Q
o30 = 4
singQ=ideal Q+ideal jacobian Q
L=saturate(singQ)
o32 = ideal (y  + 11y  + 35y , y  + y , y  - 43y  - 10y , y  + 21y  - 44y )
              1      5      6   2    6   3      5      6   4      5      6
codim(L+J),degree(L+J)
o33 = (5, 4)
\end{verbatim}
\end{small}

The $k$-normality of $X$ can be investigated by computing 
the difference between the dimension of the degree $k$ part 
of the coordinate ring {\tt T/J} of $X$ and 
$h^0(X,\cO_X(k))=-1+4k^2+3k$.
The following line will compute the Hilbert function of the coordinate 
ring of $X$ up to degree 10
(the function {\tt hilbertFunction(i,J)} returns the dimension
of the degree $i$ part of {\tt T/J} when {\tt J} is an ideal of a ring {\tt T}):
\begin{small}
\begin{verbatim}
apply(0..10,i->hilbertFunction(i,J))
o34 = (1, 6, 20, 44, 75, 114, 161, 216, 279, 350, 429)
\end{verbatim}
\end{small}
For example, for $k=1$ we see that this difference is zero,
hence $X$ is 1-normal, while for $k=2$ this difference is 1, 
hence $X$ is not 2-normal.
In this way one can check that $X$ is k-normal 
for any $k = 3,\ldots,10$.
Since it is known that any surface of the type considered in this example 
is not 2-normal, but it is $k$-normal for $k\geq 11$ (see\cite{abb1}),
the above example shows that the generic surface
of this type is in fact $k$-normal for $k\geq3$.

\begin{rem}
Given the ideal of a non-degenerate surface $X\subset\pz^r$ of degree $d$,
it follows from the Castelnuovo bound that 
$X$ is $k$-normal for $k\geq k_0=d-2+r$.
The $k$-normality for $k<k_0$ can then be checked by computing 
Hilbert function of $X$ up to degree $k_0-1$.
\end{rem}

\subsection{Second example}

Let $C$ be a smooth curve of genus $1.$ Let $E$
be a normalized rank $2$ vector bundle of degree $0,$ so we have one and
only one of the following cases:

$0)$ $E=\mathcal{O}_{C}\oplus \mathcal{O}_{C}$ and $\pz(E)=C\times \pz^{1}$

$1)$ $E=\mathcal{O}_{C}\oplus L_{0}$ where $L_{0}\neq \mathcal{O}_{C}$ but $%
\deg L_{0}=0$

$2)$ $E$ is given by the unique not trivial extension $0\rightarrow \mathcal{%
O}_{C}\rightarrow E$ $\rightarrow \mathcal{O}_{C}\rightarrow 0.$

Let us call $X_{i}$ $i=0,1,2,$ the three surfaces. It is known that if we
consider any degree $3$ divisor $B$ over $C,$ $X_{i}$ is embedded in $\pz
^{5}$ by $A=C_{0}+p^{*}B\equiv C_{0}+3f$ as a smooth scroll surface of
degree $6$ (see \cite{i1}); $g(X_{i})=1.$ In any case $C_{0}$ $\simeq C$ is
embedded as a smooth plane curve of degree $3$ and $h^{0}(X_{0},C_{0})$ $=2,$
$h^{0}(X_{i},C_{0})$ $=1$ for $i=1,2.$

About this surface we have the following proposition (see \cite{abb2})

\begin{proposition}
\label{prop2}Every $X_{i}$ is projectively normal and it is contained
exactly in only one net of quadrics $\Lambda _{i}\simeq \pz^{2}.$
Moreover: i) $\Lambda _{0}$ contains only rank $4$ quadrics whose line
vertex is generically disjoint from $X_{0}$; in $\Lambda _{0}$ there is a
smooth plane curve $\simeq C$ whose points correspond to the quadrics of $%
\Lambda _{0}$ whose vertex is contained in $X_{0}.$ ii) The generic quadric
of $\Lambda _{1}$ is smooth; the only singular quadrics in $\Lambda _{1}$
have rank $4$ and they are parametrized by a smooth plane curve $\mathcal{C}%
\simeq C;$ the discriminat divisor in $\Lambda _{1}\simeq \pz^{2}$ is a
reducible plane sextic $\mathcal{D}=2\mathcal{C}.$ iii) The generic quadric
of $\Lambda _{2}$ has rank $5$; the only rank $4$ quadrics in $\Lambda _{2}$
are parametrized by $C_{0}$: in fact their vertices are lines, tangent to $%
C_{0}$ with multiplicity $2.$
\end{proposition}

Since here $q=1$, we can take as $C$ a smooth plane cubic,
and choose $D=B$ as an effective divisors of degree $3$ 
in order to satisfy the assumptions of Theorem 
\ref{scrolls} and Remark \ref{originale}. 
As in the first example, we compute the module $H^0_*(D)$, called {\tt DS}, and
the module $H^0_*(L_0\otimes\cO_C(D))$ , called {\tt D2S}.
\begin{small}
\begin{verbatim}
K=ZZ/101;
R=K[x_0..x_2]
C=ideal random(R^1,R^{-3})
codim (C+ideal jacobian C)
Ddual=randomPoints(C,3)
D2dual=randomPoints(C,3)
S=R/C
DSdual=substitute(Ddual,S);DS=Hom(DSdual,S);
D2Sdual=substitute(D2dual,S);D2S=Hom(D2Sdual,S);
\end{verbatim}
\end{small}

\subsubsection{Case 0}
We have $M={\tt DS}\oplus {\tt DS}$. Hence we perform:
\begin{small}
\begin{verbatim}
M=DS++DS;
J=scrollIdeal(M)
dim J, degree J
o23 = (3, 8)
codim (J+minors(3,jacobian gens J)) ==6
betti res J
o25 = total: 1 7 11 6 1
          0: 1 .  . . .
          1: . 3  2 . .
          2: . 4  9 6 1
\end{verbatim}
\end{small}

By the previous set of generators it is easy to see that $X_0$ is
contained in a net of quadrics $\Lambda _{0}$.
We call {\tt Q} a set of generators for $\Lambda _{0}$:
\begin{small}
\begin{verbatim}
Q=(gens J)_{0..2}
o26 = | y_3y_5-y_2y_6 y_3y_4-y_1y_6 y_2y_4-y_1y_5 |
\end{verbatim}
\end{small}
The resolution of {\tt J} suggest that these quadrics have 2
independent linear relations among them. Indeed they are 
the $2\times2$ minors of the following matrix {\tt A},
and $X_0$ is contained in a smooth scroll of dimension 3 and degree 4 in $\pz^5$:
\begin{small}
\begin{verbatim}
A=syz Q
o27 = {2} | y_4  y_1  |
      {2} | -y_5 -y_2 |
      {2} | y_6  y_3  |
ideal Q==minors(2,A)
o28 = true
dim ideal Q, degree ideal Q
o29 = (4, 3)
codim (ideal Q+minors(codim ideal Q,jacobian Q))
o30 = 6
\end{verbatim}
\end{small}

Now we verify that $\Lambda _{0}$ contains only rank $4$ quadrics,
by checking that all quadrics have rank $\leq 4$ and no quadric
has rank $\leq 3$.
\begin{small}
\begin{verbatim}
Par=K[a,b,c];ParT=T**Par
t=substitute(vars T,ParT);par=substitute(vars Par,ParT);
genericQuadric=matrix(ParT,{{a,b,c}})*transpose substitute(Q,ParT);
matrixGenericQuadric=diff(transpose t,diff(t, genericQuadric))
matrixGenericQuadric=substitute(matrixGenericQuadric,Par)
o37 = {1} | 0  0  0 0 -c -b |
      {1} | 0  0  0 c 0  -a |
      {1} | 0  0  0 b a  0  |
      {1} | 0  c  b 0 0  0  |
      {1} | -c 0  a 0 0  0  |
      {1} | -b -a 0 0 0  0  |
Gamma5=ideal mingens minors(5,matrixGenericQuadric)
o38 = ideal 0 
Gamma4=saturate ideal mingens minors(4,matrixGenericQuadric)
o39 = ideal 1 
\end{verbatim}
\end{small}

We therefore compute the vertex locus {\tt genericVertex} of the net of quadrics in
$\pz^2\times\pz^5$ and the locus {\tt G} in $\pz^2$ of the quadrics $\Gamma_0$
whose vertex line is contained in $X_{0}$, checking that indeed
this is a smooth plane cubic:
\begin{small}
\begin{verbatim}
genericVertex=(ideal genericQuadric +ideal diff(t, genericQuadric))
W=ParT/substitute(J,ParT)
G=saturate(substitute(genericVertex,W),ideal substitute(t,W));
G=ideal mingens saturate(substitute(G,Par));
             3      2         2      3      2                 
o43 = ideal(a  + 19a b + 14a*b  - 41b  + 25a c - 44a*b*c + ...
codim ideal jacobian G
o44 = 3
\end{verbatim}
\end{small}
The fact that the curve {\tt G} is isomorphic to {\tt C} is a
geometric consequence of the construction, since each vertex line is
a line of the scroll $X_0$, which projects in a point of {\tt C}.

\subsubsection{Case 1}
Here $M$ is again a direct sum, namely ${\tt DS}\oplus {\tt D2S}$. 
\begin{small}
\begin{verbatim}
M=DS++D2S
J=scrollIdeal(M)
(dim J, degree J)
o47 = (3, 6)
codim (J+minors(3,jacobian gens J)) ==6
betti res J
o49 = total: 1 5 9 6 1
          0: 1 . . . .
          1: . 3 . . .
          2: . 2 9 6 1
\end{verbatim}
\end{small}
As in the previous subcase, we compute the representative matrix of a
generic quadric in the net $\Lambda_1$ and the the discriminant divisor,
called {\tt G2}.
Then we check that it is indeed the square of a cubic {\tt G} and that
the singular quadrics, parametrized by {\tt G}, have all rank 4:
\begin{small}
\begin{verbatim}
Q=(gens J)_{0..2}
Par=K[a,b,c];ParT=T**Par
t=substitute(vars T,ParT);par=substitute(vars Par,ParT);
genericQuadric=matrix(ParT,{{a,b,c}})*transpose substitute(Q,ParT);
matrixGenericQuadric=diff(transpose t,diff(t, genericQuadric))
matrixGenericQuadric=substitute(matrixGenericQuadric,Par)
o57 = {1} | 0            0            0            c -23a+27b+18c -47a-30b-43c |
      {1} | 0            0            0            b 21a+17b-4c   -22a-21b+14c |
      {1} | 0            0            0            a 47a+41b+10c  -12a-40b+50c |
      {1} | c            b            a            0 0            0            |
      {1} | -23a+27b+18c 21a+17b-4c   47a+41b+10c  0 0            0            |
      {1} | -47a-30b-43c -22a-21b+14c -12a-40b+50c 0 0            0            |
G2=ideal det matrixGenericQuadric
G=radical G2
codim (G+ideal jacobian G)
G==saturate minors(5,matrixGenericQuadric)
o61 = true
\end{verbatim}
\end{small}

We want to point out also the following nice geometric configuration,
not shown by Proposition \ref{prop2} and completely unexpected:
\begin{rem}
The locus $Y_1$ of the lines in $\pz^5$, which are vertices of 
the singular quadrics in the net $\Lambda_1$ of quadrics containing $X_1$,
is again a geometrically ruled surface of degree 6.
According to the classification in Proposition \ref{prop2}, 
$Y_1$ is of the same type as $X_1$.
Moreover the intersection $X_1\cap Y_1$ consists of two plane cubic curves
lying in disjoint planes, one of them being $C_0$.
\end{rem}

\begin{pf}
We aim here just to verify this configuration on the constructed example.
The locus {\tt singularVertices} 
of the lines in $\pz^5$ which are vertices of 
the singular quadrics is computable by the following commands:
\begin{small}
\begin{verbatim}
genericVertex=(ideal genericQuadric +ideal diff(t, genericQuadric))
singularVertices=saturate(genericVertex+substitute(G,ParT),ideal par);
singularVertices=saturate substitute(singularVertices,T);
dim singularVertices,degree singularVertices
o65 = (3, 6)
V=ideal mingens (singularVertices+J);
dim V,degree V
o67 = (2, 6)
\end{verbatim}
\end{small}
Therefore $Y_1$ is again a geometrically ruled surface of degree 6, 
as stated.

Now we perform the same computations as done for $X_1$.
First we verify that $Y_1$ has the same Betti numbers as $X_1$
and that the generic element of the net of quadrics 
containing $Y_1$ is smooth and that the discriminant divisor
is the square of a smooth plane cubic.
Then we verify that the intersection of $Y_1$ with $X_1$ consists 
of two plane smooth cubic curves, 
{\tt C1} and {\tt C2}.
\begin{small}
\begin{verbatim}
C1=(decompose V)_0;C2=(decompose V)_1;
betti C1,betti C2
\end{verbatim}
\end{small}
Finally, as explained in the next case, 
we compute the ideal of $C_0$ explicitly
and we verify that one of these curves is indeed $C_0$
(this computation and the check of the smoothness of $C_1$ and $C_2$
are here omitted).
\end{pf}

\subsubsection{Case 2}
Here $M$ is an extension in $Ext^1({\tt DS},{\tt DS})$, where
{\tt DS} is the module constructed as in Case 0 corresponding
to an effective divisor $D$ of degree 3:  
\begin{small}
\begin{verbatim}
M=randomExt(DS,DS)
\end{verbatim}
\end{small}
where {\tt randomExt()} is the function defined in section 3.1.

Unfortunately, this function returns an error, revealing a 
not correctly defined code for the function {\tt random(Module,Module)}, 
contained in {\tt Macaulay2} system.
We therefore have to correct the definition of the function
{\tt random(Module,Module)}, which appears
in the file {\tt genmat.m2} of {\tt Macaulay2} package:
copy its definition and the one of {\tt randommat()}
into a file, replace the last two lines of
\begin{small}
\begin{verbatim}
          else (
               m := basis(deg,R);
               s := degreesTally#deg;
               reshape(F,G, 
                    m * randommat(R, numgens source m, s))))
\end{verbatim}
\end{small}
with the correction
\begin{small}
\begin{verbatim}
               map(F,G, reshape(cover F,G, 
                    m * randommat(R, numgens source m, s)))))
\end{verbatim}
\end{small}
and redefine both functions 
{\tt randommat()} and {\tt random(Module,Module)} in
the program {\tt Macaulay2}.
Another possibility is to replace all the 5 lines quoted above
with a line with only the character {\tt )}:
without this {\tt else} subcase the function {\tt random(Module,Module)} works correctly.
Once this correction is done, the command
{\tt M=randomExt(DS,DS)} will work.


We proceed now to compute the ideal {\tt J} of $X_1\subset \pz^5$,
to check the smoothness and to give the Betti table of $I_{X_1}$:
\begin{small}
\begin{verbatim}
J=scrollIdeal(M)
dim J, degree J
o73 = (3, 6)
codim (J+minors(3,jacobian gens J)) ==6
betti res J
o75 = total: 1 5 9 6 1
          0: 1 . . . .
          1: . 3 . . .
          2: . 2 9 6 1
\end{verbatim}
\end{small}

As in the previous subcase, we compute the representative matrix of a
generic quadric in the net $\Lambda_2$, we check that all the
quadrics have rank $\leq 5$, and we compute the divisor of
the rank 4 quadrics in $\Lambda_2$, a smooth cubic {\tt G}:
\begin{small}
\begin{verbatim}
Q=(gens J)_{0..2}
Par=K[a,b,c];ParT=T**Par
t=substitute(vars T,ParT);par=substitute(vars Par,ParT);
genericQuadric=matrix(ParT,{{a,b,c}})*transpose substitute(Q,ParT);
matrixGenericQuadric=diff(transpose t,diff(t, genericQuadric))
matrixGenericQuadric=substitute(matrixGenericQuadric,Par)
o83 = {1} | 2c           b           -45a-18b+31c 0           -34a-37b-35c -23a+30b+44c |
      {1} | b            2a          33a-5b-c     34a+37b+35c 0            -3a+14b+c    |
      {1} | -45a-18b+31c 33a-5b-c    -49a-40b-19c 23a-30b-44c 3a-14b-c     0            |
      {1} | 0            34a+37b+35c 23a-30b-44c  0           0            0            |
      {1} | -34a-37b-35c 0           3a-14b-c     0           0            0            |
      {1} | -23a+30b+44c -3a+14b+c   0            0           0            0            |
det matrixGenericQuadric
o84 = 0
G=saturate minors(5,matrixGenericQuadric)
             3      2         2      3      2                
o85 = ideal(a  - 14a b + 44a*b  - 32b  - 48a c - 30a*b*c + ...
codim ideal jacobian G
\end{verbatim}
\end{small}

\medskip
We want now to verify, according to Proposition \ref{prop2}, that the
vertex of any rank 4 quadric in $\Gamma_2$ is a line tangent to
$X_2$ at a point of $C_0$. Again, we obtain further results similar
to the previous case:

\begin{rem}
The locus $Y_2$ of the lines in $\pz^5$, which are vertices of 
the singular quadrics in the net $\Lambda_2$ of quadrics containing $X_2$,
is again a geometrically ruled surface of degree 6.
According to the classification in \ref{prop2}, 
$Y_2$ is of the same type as $X_2$.
Moreover the intersection $X_2\cap Y_2$ consists of the
cubic $C_0$ counted twice.
Therefore a vertex of a quadric in $\Lambda_2$ is a line $L$ tangent to
$X_2$ in the point of $C_0$ given by the intersection of $L$ with $C_0$.
\end{rem}

\begin{pf}
We aim here just to verify this configuration on the constructed example.
First we give a function to compute the ideal of the fiber in $\pz^5$ of an
effective divisor over $C$ and we compute the ideal of $C_0$, called {\tt C0}:
\begin{small}
\begin{verbatim}
pullbackIdeal = (I) -> (
     R:=ring I;TR:=ring IS;        
     J:=substitute(I,TR)+IS;
     J=saturate(J,ideal substitute (vars R, TR));
     ideal mingens substitute(J,T))
betti (H=pullbackIdeal(Ddual))
dim H, degree H
C0=(ideal H_0+J):H
\end{verbatim}
\end{small}

Then we compute the locus $Y_2$.
we intersect this surface with $X_2$, 
and we verify that the intersection $Y_2\cap X_2$ is
the square of $C_0$:
\begin{small}
\begin{verbatim}
genericVertex=(ideal genericQuadric +ideal diff(t, genericQuadric))
singularVertices=saturate(genericVertex+substitute(G,ParT),ideal par);
singularVertices=saturate substitute(singularVertices,T);
dim singularVertices,degree singularVertices
o93 = (3, 6)
V=ideal mingens (singularVertices+J);
dim V,degree V
o95 = (2, 6)
betti (C1=radical V)
o96 = generators: total: 1 4
                      0: 1 3
                      1: . .
                      2: . 1
C1==C0
o97 = true
\end{verbatim}
\end{small}
The check that $Y_2$ is again a geometrically ruled surface of degree 6
and that the generic element of the net of quadrics 
containing $Y_2$ has rank 5 is similar to the one already 
performed for $X_2$.

The last statement is a direct consequence of the previous ones:
any line $L$, vertex of a rank four quadric of $\Lambda_2$, 
intersects $X_2$ at two points belonging to $C_0$. 
Since $L$ cannot lie in the plane containing $C_0$ (otherwise $L$
would cut $C_0$, hence also $X_2$, in a dimension 0 scheme of degree
3, contrasting with the above calculation for {\tt V}), 
it follows immediately that $L$ is a tangent line to $X_2$ at
a unique point of $C_0$, the intersection of $L$ with $C_0$.
\end{pf}

\section{Conic bundles}

In this section we consider the explicit construction of surfaces which are $%
\pz^{1}$-bundles over a smooth curve $C,$ embedded in such a way that
every fibre is a smooth conic. In other words, we consider 
polarized surfaces $(X,A)$ as in \S 1, 
where $X:=\pz(E)$ and $A=2C_{0}+p^{*}B$ is very ample. 
The aim is therefore to prove the following:

\begin{bigthm}
\label{conic bundles} Let $C\subset \pz^m$ be a smooth curve $C$ of genus $q$,
$B$ a divisor on $C$ and $L$ a line bundle over $C$.
Consider a normalized rank 2 vector bundle $E\in\mathcal Ext^1(L,\cO_C)$ 
over $C$ given by an extension $0\ra \cO_C \ra E \ra L \ra 0$ and
suppose that the divisor $A=2C_{0}+p^{*}B $ on the surface
$X=\pz(E)$ is very ample. 
Then there is an algorithm yielding a set of generators for the ideal $I_{X}$ 
of the embedded $X$ in $\pz^{h^{0}(X,A)-1}=\pz(H^{0}(X,A)^{*})$ by $|A|$.
\end{bigthm}

Again, before giving the proof of theorem, in terms of an explicit
algorithm, let us fist develop a technical criterium needed in the
given algorithm.

\medskip
We firstly cite the following theorem of Butler:

\begin{theorem}\label{5.1}
\cite[Thm. 5.1A]{bu} $(\rm{char}\  \kz=0$ or $q\leq 1).$ 
Let $E$ be a vector bundle over a smooth projective curve $C$ of genus
$q$, $p: E \ra C$ the projection, and let $X=\pz(E)$. 
If $Z$ is a $(-1)$ $p$-regular line bundle over $X$ with
$\mu^{-}(p_*Z)>2q$,
then $Z$ is normally generated.
\end{theorem}

We will use this result in the form of the following corollary:

\begin{cor} \label{cor n-norm}
$(\rm{char}\  \kz=0).$
Let $E$ be a vector bundle over a smooth projective curve $C$ of genus
$q$, $p: E \ra C$ the projection, and let $X=\pz(E)$. 
If $\mu^{-}(E)>2q$,
the tautological divisor $\tau$ of $X$ is very ample and
$X\subset \pz(H^0(X,\tau)^*)$ is projectively normal.
\end{cor}

\begin{pf}
It is well known that the condition $\mu^{-}(E)>2q$
implies that the tautological divisor $\tau$ of $X$ is very ample,
cf. Lemma 1.12 of \cite{bu}.
Since $\tau$ is very ample, it is enough to prove that $X$ is normally
generated in $\pz(H^0(X,\tau)^*)$.
For this we apply Theorem \ref{5.1}.

We recall that a divisor $Z$ is called $(-1)$ $p$-regular if,
for every fibre $f$ of $p$ over $C$, $H^i(f,Z_{|f}(-1-i))=0$,
for all $i>0$.
In our case these groups are $H^i(\pz^r,\cO_{\pz^r}(-i))=0$, where
$r=\rk E$,
hence $\tau$ is automatically $(-1)$ $p$-regular.

For the second condition, we have $\mu^{-}(p_*\tau)=\mu^{-}(E)$.
\end{pf}

We are now ready to give the required criterium.

\begin{proposition}\label{X' norm}
  Let $E$ be a vector bundle over a (smooth) curve 
$C$ of genus $q$ and $D$ an effective divisor of degree $d$ on $C$.
If the condition
\begin{equation}\label{deg2norm}
 \mu^{-}(E)+d>2q
\end{equation}
is satisfied, the divisor $C_0+p^*(D)$ is very ample on $X=\pz(E)$ 
and the image $X'$ of $X$, given by the linear system
$|C_0+p^*D|$, is projectively normal.
\end{proposition}

\begin{pf}
  Just apply Corollary \ref{cor n-norm} to $E':=E\otimes\cO_C(D)$
and recall that $\mu^{-}(E')=\mu^{-}(E)+d$.
\end{pf}

\begin{rem}
Condition (\ref{deg2norm}) is not an evident numerical condition,
since it is not clear how to compute $\mu^{-}(E)$ for a given
vector bundle $E$.
 
However, if the genus $q$ of the curve $C$ satisfies $q\geq 2$, 
the set of points in $Ext^1(L,\cO_C)$ parametrizing a semi-stable
vector bundle $E$  is a Zariski open set, see the classical \cite[Thm. 2]{ns}.
Hence for a general choice of such an extension the corresponding $E$ is 
semistable and $\mu^-(E)=\mu(E)=-\frac{e}2$.

For the case $q=1$ it is known that if $E$ is indecomposable, then $E$ is semi-stable.
If instead $E$ is decomposable, say $E=L\oplus L'$, then 
$\mu^-(E)=\min(\deg(L),\deg(L'))$, see \cite[Lemma 2.8]{abb3}.
In the case $q=0$, $E$ is necessarily of the type $E=\cO(a)\oplus\cO(b)$ and
$\mu^-(E)=\min(a,b)$.
\end{rem}

\medskip
\begin{pf}[Proof of Thm. \ref{conic bundles}]
Chose an effective divisor $D$ on $C$ of degree $d$ such that $D-B$ is
effective and such that $D$ satisfies the condition (\ref{deg2norm}).
Then the divisor $C_0+p^*(D)$ on the surface $X=\pz(E)$ is very ample 
and, letting $X'\subset\pz^r$ be the image of the embedding 
$\iota:\ X=\pz(E) \hookrightarrow \pz^r$ 
given by $|C_0+p^*(D)|$, 
the surface $X'$ is projectively normal by Proposition \ref{X' norm}.
By applying Theorem \ref{scrolls},
we obtain a set of generators for the ideal $I_{X'}$ of $X'$ in $\pz^r$.

Let $R$ be a section of the sheaf $i_*p^*\cO_C(2D-B)$ and
let $H$ be the hyperplane divisor of $\pz^r$.
Notice that $R$ is an effective divisor of $X'$ and that 
we have
$$0 \ra \cI_{X'} \ra \cI_{R} \ra \cI_{R,\,X'} \ra 0,$$
where $\cI_{R,\,X'}$ is the relative ideal sheaf of $R$ in $X'$.

Recall that $X'$ is 2-normal.
Therefore $H^1(\pz^r,\cI_{X'}(2H))=0$ and
the above sequence, tensorized with $\cO_{\pz^r}(2H)$, gives
$$0 \ra H^0(\pz^r,\cI_{X'}(2H)) 
\ra H^0(\pz^r,\cI_{R}(2H))
\ra H^0(\pz^r,\cI_{R,\,X'}(2H)) \ra 0.$$
Since 
$\cI_{R,\,X'}(2H)\cong \cO_{X'}(2H-R)\cong \cO_X(2C_0+p^*(B))=\cO_X(A)$, 
we have 
$$H^0(X,A)\cong \frac{H^0(\pz^r,\cI_{R}(2H))}{H^0(\pz^r,\cI_{X'}(2H))}.$$

Let $f_0,\dots,f_n$ be a set of representative quadrics in 
$H^0(\pz^r,\cO_{\pz^r}(2))$
for a basis of the above quotient space, where $n=h^0(X,A)-1$.
The image of $X$ under the linear system $|A|$ is then given
in the following way:
if $y_0,\ldots y_n$ are indeterminates and $S$ is
the coordinate ring of $X'\subset\pz^r$,
the ideal $I_X$ is the kernel of the map
$K[y_0,\ldots,y_n] \ra S$ obtained by sending $y_i$ to
$[f_i]$, where $[f_i]$ is the class of $f_i$ in $S$.
\end{pf}

\medskip

\subsection{\bf The 3-fold $\pz(S^2E')$, the surface 
$\pz(E')$ and its 2-Veronese image.}

In this subsection we want to describe shortly some varieties related
to $X'$, as the title suggests, and to describe their geometric correlations.
According to the previous notation, 
let $D$ satisfy the hypothesis (\ref{deg2norm}) and $E'=E\otimes\cO_C(D)$.

\medskip
Let us consider the 2-Veronese embedding 
$\nu:\ \pz^r \hookrightarrow \pz^N$,
where $N=\binom{r+2}2-1$
and let $X''$ be the image of $X'$ under $\nu$,
i.e. the image of $X$ under the composition $\nu\circ\iota$.
Then $X''$ is the image of $X$ via the map associated to the 
linear system $|2C_0+p^*(2D)|$.
Algebraically, the map is given as follows.
Let $y_0,\dots,y_{r}$ be a 
base of $H^0(C,E')$: 
then $y_0^2,y_0y_1,\dots,y_{r}^2$ is a base
of $S^2(H^0(C,E'))$ and,
denoting $z_{i,j}$ with $i\leq j$ as a set of coordinates for $\pz^N$,
the composition $\nu\circ\iota$ is 
given by mapping $z_{i,j}$ to the product $y_i y_j$,
considered as an element in $\cO_{\pz(E')}(2)$.

Moreover, the map $\nu^*$ gives an isomorphism
between the vector space of  hyperplanes 
containing the image of the effective divisor $R$, 
and the vector space $H^0(\pz^r,\cI_{R}(2H))$,
where $R$ is the divisor used in the proof of Theorem 
\ref{conic bundles}, i.e. a section of $i_*p^*(2D-B)$.

\medskip
Define now $E_1:=S^2(E')$ and consider
the 3-fold $\pz(E_1)$:
we want to compute the ideal of the image $X_1$ of $\pz(E_1)$ in
$\pz^s:=\pz(H^0(C,E_1)^*)\subset \pz^N$ via the linear system
given by the tautological divisor $T_1$ of $E_1$.

From our hypothesis (\ref{deg2norm}) over $D$,
it is easy to see that also $S^2(E')$ is very ample,
e.g. because $\mu^{-}(S^2(E'))=2\mu^{-}(E')>4q$ and Lemma 1.12 of \cite{bu}.
Since $X'$ is projectively normal, 
$S^2(H^0(C,E'))$ surjects to $H^0(C,S^2(E'))$,
the kernel being $H^0(\pz^r,\cI_{X'}(2))$.
Indeed, from the projective normality we obtain the exact sequence 
\begin{equation}\label{pn}
0 \ra H^0(\pz^r,\cI_{X'}(2)) \ra H^0(\pz^r,\cO_{\pz^r}(2)) \ra 
H^0(X',\cO_{X'}(2)) \ra 0,
\end{equation}
where $H^0(\pz^r,\cO_{\pz^r}(2))\cong S^2(H^0(C,E'))$ and 
$H^0(X',\cO_{X'}(2))\cong H^0(C,S^2(E'))$.

For simplicity, let us now assume that $C\subset \pz^m$ satisfies
\begin{equation}\label{S-gen}
\deg(C)\geq 2q+1,
\end{equation}
so that $C$ is projectively normal.
Then $H^0_*(C,E')$ is generated by $H^0(C,E')$, as module over the ring
$\oplus_{t\geq 0}H^0(C,\cO_C(t))$.
Indeed by the above hypothesis we have both $\deg (\cO_C(t))\geq 2q$ for
$t\geq 1$ and $\mu^{-}(E')>2q$, and the surjectivity of the natural map
$H^0(C,E')\otimes H^0(C,\cO_C(t)) \ra H^0(C,E'(t))$ is an application
of Theorem 2.1 of \cite{bu}.
In the same way, $H^0(C,S^2(E'))$ generates the module $H^0_*(C,S^2(E'))$.

\smallskip
We claim that $S^2(H^0_*(C,E'))$ surjects to $H^0_*(C,S^2(E'))$.
Indeed, if
\begin{equation}
  \label{h0}
M_2\ra M_1 \ra H^0_*(C,E')\ra 0 
\end{equation}
is a free presentation over $S$ for the module $H^0_*(C,E')$
with $\rk M_1=h^0(C,E')$,
then the following is a free presentation for $S^2(H^0_*(C,E'))$ with 
$\rk S^2(M_1)=\binom{h^0(C,E')}2$:
\begin{equation}
  \label{s2h0}
(M_2\otimes M_1) \rxa{\phi} S^2 M_1 \ra S^2(H^0_*(C,E')) \ra 0.
  \end{equation}
Moreover, because of (\ref{S-gen}), 
$y_0,\dots,y_r$ is a set of generators of the module $H^0_*(C,E')$ and therefore
$y_0^2,y_0y_1,\dots,y_r^2$ is a set of generators of $S^2(H^0_*(C,E'))$.
By the projective normality, 
the images of $y_0^2,y_0y_1,\dots,y_r^2$
in $H^0_*(C,S^2(E'))$ span the whole vector space $H^0(C,S^2(E'))$, 
which generates $H^0_*(C,S^2(E'))$, again by (\ref{S-gen}),
and the claim is proved.

\smallskip
Now we compute the kernel of the natural map
$S^2(H^0_*(C,E'))\ra H^0_*(C,S^2(E'))$.
Let $N$ be the set of the linear forms $\sum\alpha_{i,j} z_{i,j}$ in 
$S^2(H^0(C,E'))$ such that $\sum \alpha_{i,j} y_i y_j \in H^0(\pz^r,I_{X'}(2))$.
By (\ref{pn}), the symmetric algebra $S(H^0(C,S^2(E')))$ over the 
vector space $H^0(C, S^2(E'))$ is isomorphic to the symmetric
algebra  $S(V)$ over $V:=S^2(H^0(C,E'))/N$.
Hence a polynomial $f\in S(S^2(H^0(C,E')))$ is zero in
$S(H^0(C,S^2(E')))$ if and only if $f$ is in the ideal generated by
$N$.
Thus the kernel of the natural map from $S(S^2(H^0_*(C,E')))$ to
$S(H^0_*(C,S^2(E')))$ defined by $z_{i,j}\mapsto y_iy_j$
is exactly the ideal generated by $N$ in $\kz[z_{i,j}]$
(cf. assumption (\ref{S-gen})).

\medskip
We are therefore able, as in the proof of Lemma \ref{risoluzione},
to get a set of generators for the ideal of $X_1$ in $\pz^N$:
we multiply again the map $\phi$ in (\ref{s2h0}) with the variables
$z_{i,j}$ (by considering their exact order), and we 
also add the linear forms in $N$ to these equations.
Then we saturate as usual with respect to the irrelevant 
ideal $(x_0,\ldots,x_m)$ of $\pz^m$.
Notice that $X_1$ is degenerate if $H^0(\pz^r,I_{X'}(2))\neq 0$: 
$X_1$ lies in the $\pz^s\subset \pz^N$ given by
equations in $N$, obtained from $H^0(\pz^r,I_{X'}(2))$.

\begin{rem}
If $\deg(C)\leq 2q$ and $C\subset \pz^m$ is not projectively normal, 
then we can still find explicit equations of $X_1$,
by arguing as already done in the proof of Lemma \ref{risoluzione}. 
Indeed, instead of the presentation (\ref{h0}), 
we take a presentation for the submodule $M'\subset H^0_*(C,E')$ 
generated by $H^0(C,E')$ with $\rk M_1=h^0(C,E')$,
and the corresponding presentation $\phi$ of the submodule
generated by $S^2(H^0(C,E'))$ in $S^2(H^0_*(C,E'))$.
Since $X'$ is projectively normal,  
$S^2(H^0(C,E'))$ surjects to $H^0(C,S^2(E'))$
and the rest of the algorithm works without further modifications.
\end{rem}

\medskip
The ideal of the surface $X''$ is obtained by the ideal
of $X_1$ by adding to it the quadrics of $\pz^N$
in the kernel of the map
$$S^2(S^2(H^0(C,E')))\ra H^0(C,S^4(E')).$$
These quadrics indeed generate all relations among the $y_iy_j$,
the images of the $z_{i,j}$ in the $S$-algebra $S(H^0_*(C,E'))$.
Considering sheaves, we have
$0 \ra \cL \ra S^2(S^2(E')) \ra S^4(E') \ra 0$:
$\cL$ is a line bundle ($E'$ has rank 2)
and the quadrics obtained above are global sections of $\cL$.
We call these quadrics the {\em relative Veronese quadrics}, since
they give fiberwise the quadric ideal of the Veronese embedding 
$\pz^1 \ra \pz^2$ or zero.

\subsection{A linearly embedded surface and an alternative proof of
  Thm. \ref{conic bundles}}

From the short exact sequence
$0\ra \cO_C \ra E \ra L \ra 0$
one can derive two other short exact sequences where $E_1$ sits,
namely:
\begin{equation}\label{eqE_1 a}
  0\ra E\otimes \cO_C(2D) \ra E_1 \ra L^2\otimes\cO_C(2D) \ra 0
\end{equation}
and 
\begin{equation}\label{eqE_1 b} 
  0 \ra \cO_C(2D) \ra E_1 \ra E\otimes L\otimes \cO_C(2D) \ra 0.
\end{equation}

Indeed, we have $0\ra E\otimes\cO_C \ra S^{2}(E) \ra S^2(L)\ra 0$,
which
can be rewritten as 
$0\ra E \ra S^{2}(E) \ra L^2 \ra 0$.
From this sequence we can proceed in two ways:
either we tensorize with $\cO_C(2D)$, getting the first claimed
sequence, or we dualize it, getting
$0\ra L^{-2} \ra S^2(E^*) \ra E^* \ra 0$,
where $L^{-i}$ denotes the $i$-th tensor power of $L^*$.
Since $E^*\cong E\otimes L^*$, we obtain the second claimed sequence 
tensorizing 
with $L^{2}\otimes\cO_C(2D)$.

\medskip
We want to use here the second sequence.
If $\mu^-(E)+2d>2q-e$, then, by Lemma 1.12 of \cite{bu},
$E\otimes L\otimes \cO_C(2D)$ is very ample and its tautological
bundle $\tau$ defines an embedding $\Sigma$ of $\pz(E)$.

Note that $T_1\mid_\Sigma=\tau$, and the ideal of $\Sigma$ in $X_1$ 
is generated by the elements in $H^0(X_1,\cO_{X_1}(T_1-{p_1}^*(2D)))$, 
where $p_1$ is the projection of $X_1$ to $C$. 

We assume furthermore $d>q-1$:
then $\Sigma$ is linearly normal in $X_1$, because $h^1(C,2D)=0$.

We can also proceed with $\Sigma$ to compute a set of generators for the ideal
$I_X$ of the image of $X$ via the map associated to the linear system
$|A|$
of Theorem \ref{conic bundles}, 
but then further conditions on $D$ are needed to be satisfied.

Suppose that $D$ satisfy, besides the (already required) conditions
$\mu^-(E)+d>2q$, $\mu^-(E)+2d>2q-e$ and $d>q-1$,
the further condition
\begin{equation}
  \label{d++1}
     H^0(C,L^{2}\otimes\cO_C(4D-B))\neq 0,
\end{equation}
where $B$ is as in the statement of Theorem \ref{conic bundles}.
Choosen a divisor 
$D'$ corresponding to $L\otimes\cO_C(2D)$, we
can write $\tau=C_0+p^*(D')$.
Then the sheaf $\cO_{X_1}(2T_1-p_1^*(2D'-B))$ on $X_1$
restricts to $\Sigma$ to the sheaf
$$\cO_\Sigma(2\tau-p^*(2D'-B))\cong
\cO_{X}(2(C_0+p^*D')-p^*(2D'-B))\cong
\cO_{X}(2C_0+p^*B)=\cO_X(A).$$

$X_1$ is projectively normal, by applying again Corollary
\ref{cor n-norm}, since  $\mu^-(E_1)=2\mu^-(E)+2d>4q$.
Hence, if
\begin{equation} \label{d++2}
   |2T_1-p_1^*(2D'-B)|\ra |2\tau-p^*(2D'-B)|
\text{ is surjective},
\end{equation}
we can proceed analogously to the proof of Theorem \ref{conic bundles}: 
considered an effective divisor $R_1\in H^0(X_1,p_1^*(2D'-B))$,
we have 
$$H^0(X,\cO_{X}(A))\cong \frac{H^0(\pz^s,\cI_{R_1}(2H_1))}
{H^0(\Sigma,\cI_\Sigma(2H_1))},$$ 
where $H_1$ denotes an hyperplane of $\pz^s\subset \pz^N$,
and again it is straightforward to compute $I_X$.

\begin{rem}
The restriction map $|2T_1-p_1^*(2D'-B)|\ra |2\tau-p^*(2D'-B)|$
is surjective if one of the following conditions are satisfied:
\begin{equation}
  \begin{cases}
    h^1(\cO_C(B))=h^1(\cO_C(D+B)\otimes L^{-2})=h^1(\cO_C(D+B)\otimes
    L^{-1})=0\\
    h^1(\cO_C(B)\otimes L^{-2})=h^1(\cO_C(B)\otimes L^{-1})=h^1(\cO_C(B))=0\\
    2\mu^-(E)+\deg(B)>2q-2+2e
  \end{cases},
\end{equation}
where $L^{-i}$ denotes the $i$-th tensor power of $L^*$.
\end{rem}

\begin{pf}
Recall that $\Sigma\in\mid T_1-p_1^*(\cO_C(2D))\mid$. 
Therefore 
$H^1(X_1,\cI_\Sigma(2T_1-p_1^*(2D'-B))=
H^1(X_1,\cO_\Sigma(T_1-p_1^*(\cO_C(2D-B)\otimes L^2)))\cong 
H^1(C,E_1\otimes \cO_C(B-2D)\otimes L^{-2})$ and the surjectivity follows from
the vanishing of $H^1(C,E_1\otimes \cO_C(B-2D)\otimes L^{-2})$.

Tensorizing the exact sequences (\ref{eqE_1 a}) and (\ref{eqE_1 b}) with 
$\cO_C(B-2D)\otimes L^{-2}$ and using appropriate tensorizations
of the sequence $0\ra \cO_C \ra E \ra L \ra 0$ we get the fisrt two  ways of
obtaining the desired vanishing.

Alternatively, one can again work with $\mu^-(E)$.
By Lemma 1.12 of \cite{bu}, it is enough that 
$\mu^-(E_1\otimes \cO_C(B-2D)\otimes L^{-2})>2q-2$, 
i.e. the third condition.
\end{pf}

\section{An Example of  conic bundles}
Let $C$ be a smooth curve of genus $1.$ 
Let $E$ be a normalized rank $2$ vector bundle of 
degree $1$ which is given by the only non trivial extension:

\begin{center}
$0\ra \cO_C \ra E \ra \cO_C(P) \ra 0$
\end{center}
where $\det (E)=c_1(E)=\cO_C(P)$ and $P$ is a fixed point of $C$.

Let $Q$ be any other point of $C$, eventually $Q=P$.
It is known that, on the surface $X=\pz(E)$, 
the divisor $A=2C_{0}+p^{*}Q\equiv 2C_{0}+f$ is very ample, 
whatever $Q$ is chosen, 
$h^{0}(X,A)=h^{0}(C,S^{2}(E)\otimes \cO_C(P))=6$ and $A$ 
embeds $X$ in $\pz^{5}$ as a smooth ruled surface of degree $8$, 
whose fibres are embedded as smooth plane conics 
(see \cite{i1}); $g(X)=3.$ In this case, $C_{0}$ is embedded as a smooth
plane cubic and $h^{0}(X,C_{0})$ $=1.$

\medskip
According to Theorem \ref{conic bundles}, we have $B=Q$ and we can choose
$D=Q+Q'$, where $Q'$ is a further point, 
so that $D-B$ is effective. The divisor $2D-B$ will be $Q+2Q'$.

\medskip
As usual we perform the necessary steps to obtain
$M=H^0_*(C,E\otimes\cO_C(D))$:
we fix the choice of the curve $C$, named {\tt C}, and 
of the three points $P,Q,Q'$, named respectively 
{\tt p}, {\tt q}, {\tt q'}.
\begin{small}
\begin{verbatim}
K=QQ;
R=K[x_0..x_2]
C=ideal (x_0*(x_2)^2-x_1*(x_1+x_0)*(x_1+2*x_0))
p=ideal (x_1,x_2);q=ideal (x_1,x_0);q'=ideal (x_1+x_0,x_2);
--C=ideal random(R^1,R^{-3})
--p=randomPoints(C,1),q=randomPoints(C,1),q'=randomPoints(C,1)
Ddual=intersect(q,q')
D2dual=intersect(Pdual,Ddual)

S=R/C
DSdual=substitute(Ddual,S);DS=Hom(DSdual,S);
D2Sdual=substitute(D2dual,S);D2S=Hom(D2Sdual,S);
M=randomExt(D2S,DS)
\end{verbatim}
\end{small}

Now we compute the ideal {\tt J} of 
$X'\subset\pz^4$,
the embedding of $\pz(E)$ through the linear system $|C_0+p^*(D)|$:
\begin{small}
\begin{verbatim}
J=scrollIdeal(M) --ideal of X'
dim J, degree J
o25 = (3, 5)
codim (J+minors(2,jacobian gens J)) ==5
betti res J
o27 = total: 1 5 5 1
          0: 1 . . .
          1: . . . .
          2: . 5 5 1
\end{verbatim}
\end{small}

The vector space $H^0(\pz^4,\cI_R(2H))$ is therefore given 
by the quadrics of $\pz^4$ passing through the fibers over $Q+2Q'$.
Moreover, since $I_{X'}$ has no elements of degree 2,
this vector space gives exactly $H^0(X,A)$.
\begin{small}
\begin{verbatim}
betti (q'squareFiber=pullbackIdeal(q'^2))
betti (qFiber=pullbackIdeal(q))
betti (A=intersect(q'squareFiber,qFiber))
o30 = generators: total: 1 7
                      0: 1 .
                      1: . 6
                      2: . 1
Q=super basis(2,A) --the linear system |2H-2D+B|
\end{verbatim}
\end{small}

At this point, it remains only to compute the image of $X'$
via the embedding given by $H^0(X,A)$, which is standard.
We report here also the Betti table and the Hilbert function of $X$.
\begin{small}
\begin{verbatim}
Z=K[z_0..z_5];S'=T/J;f=map(S',Z,substitute(Q,S'))
I=ker f --ideal of X
dim I, degree I
o36 = (3, 8)
codim (I+minors(3,jacobian gens I)) ==6
o37 = true
betti res I
o38 = total: 1 9 15 8 1
          0: 1 .  . . .
          1: . 1  . . .
          2: . 8 15 8 1
apply(0..10,i->hilbertFunction(i,I))
o39 = (1, 6, 20, 42, 72, 110, 156, 210, 272, 342, 420)
\end{verbatim}
\end{small}

\section{Embeddings with fibers of higher degree}

In this section we consider the explicit construction of surfaces
which are $\pz^{1}$-bundles over a smooth curve $C,$ embedded in such a way that
every fibre is a rational curve of degree $k\geq 3$. 
In other words, we consider 
polarized surfaces $(X,A)$ as in \S 1, 
where $X:=\pz(E)$ and $A=kC_{0}+p^{*}B$ is very ample. 
The aim is therefore to prove the following:

\begin{bigthm}
\label{k-bundles} Let $C\subset \pz^m$ be a smooth curve $C$ of genus $q$,
$B$ a divisor on $C$ and $L$ a line bundle over $C$.
Consider a normalized rank 2 vector bundle $E\in\mathcal Ext^1(L,\cO_C)$ 
over $C$ given by an extension $0\ra \cO_C \ra E \ra L \ra 0$ and
suppose that the divisor $A=kC_{0}+p^{*}B $ on the surface
$X=\pz(E)$ is very ample. 
Then there is an algorithm yielding a set of generators for the ideal $I_{X}$ 
of the embedded $X$ in $\pz^{h^{0}(X,A)-1}=\pz(H^{0}(X,A)^{*})$ by $|A|$.
\end{bigthm}

\begin{pf}
Analogous to the one of Theorem \ref{conic bundles}.
Chose an effective divisor $D$ on $C$ of degree $d$ such that $D-B$ is
effective and such that $D$ satisfies the condition (\ref{deg2norm}).
Then the divisor $C_0+p^*(D)$ on the surface $X=\pz(E)$ is very ample 
and, letting $X'\subset\pz^r$ be the image of the embedding 
$\iota:\ X=\pz(E) \hookrightarrow \pz^r$ 
given by $|C_0+p^*(D)|$, 
the surface $X'$ is projectively normal by Proposition \ref{X' norm}.
By applying Theorem \ref{scrolls},
we obtain a set of generators for the ideal $I_{X'}$ of $X'$ in $\pz^r$.

Let $R$ be a section of the sheaf $i_*p^*\cO_C(kD-B)$ and
let $H$ be the hyperplane divisor of $\pz^r$.
Notice that $R$ is an effective divisor of $X'$ and that 
we have
$$0 \ra \cI_{X'} \ra \cI_{R} \ra \cI_{R,\,X'} \ra 0,$$
where $\cI_{R,\,X'}$ is the relative ideal sheaf of $R$ in $X'$.

Recall that $X'$ is $k$-normal.
Therefore $H^1(\pz^r,\cI_{X'}(kH))=0$ and
the above sequence, tensorized with $\cO_{\pz^r}(kH)$, gives
$$0 \ra H^0(\pz^r,\cI_{X'}(kH)) 
\ra H^0(\pz^r,\cI_{R}(kH))
\ra H^0(\pz^r,\cI_{R,\,X'}(kH)) \ra 0.$$
Since 
$\cI_{R,\,X'}(kH)\cong \cO_{X'}(kH-R)\cong \cO_X(kC_0+p^*(B))=\cO_X(A)$, 
we have 
$$H^0(X,A)\cong \frac{H^0(\pz^r,\cI_{R}(kH))}{H^0(\pz^r,\cI_{X'}(kH))}.$$

Let $f_0,\dots,f_n$ be a set of representative degree $k$ forms in 
$H^0(\pz^r,\cO_{\pz^r}(k))$
for a basis of the above quotient space, where $n=h^0(X,A)-1$.
The image of $X$ under the linear system $|A|$ is then given
in the following way:
if $y_0,\ldots y_n$ are indeterminates and $S$ is
the coordinate ring of $X'\subset\pz^r$,
the ideal $I_X$ is the kernel of the map
$K[y_0,\ldots,y_n] \ra S$ obtained by sending $y_i$ to
$[f_i]$, where $[f_i]$ is the class of $f_i$ in $S$.
\end{pf}

\begin{rem}
A difficult point is to compute the system of hypersurfaces of degree $k$ in
$\pz^r$ through $\deg(kD-B)$ lines of the scroll $X'$.
We do not know how hard is this task computationally
when $k$ or the degree of $D$ increases (increasing the degree of $D$
increases also $r$ and the complexity of the computation). 
\end{rem}

\section{An example of surfaces with higher degree fibers}

Let $C$ be a smooth curve of genus $1.$ 
Let $E$ be a normalized rank $2$ vector bundle of 
degree $1$ which is given by the only non trivial extension:
\begin{center}
$0\ra \cO_C \ra E \ra \cO_C(P) \ra 0$
\end{center}
where $\det (E)=c_1(E)=\cO_C(P)$ and $P$ is a fixed point of $C$.

On the surface $X=\pz(E)$, the divisor $A=3C_{0}$ is very ample:
the numerical criterion of Reider is satisfied (cf. \cite{rei}).
$h^{0}(X,A)=h^{0}(C,S^{3}(E))=6$ and $A$ 
embeds $X$ in $\pz^{5}$ as a smooth ruled surface of degree $9$, 
whose fibres are embedded as twisted cubics (see \cite{i1}); 
$g(X)=4$. 
In this case, $C_{0}$ is embedded as a smooth
plane cubic as $h^{0}(X,2C_{0})=3.$

\smallskip
According to Theorem \ref{k-bundles}, we have $B=0$ and
we can choose $D=Q+Q'$, where $Q,Q'$ are any pair of points.
The divisor $3D-B$ will then be $3D=3Q+3Q'$.

\smallskip
As in section 4, we perform the necessary steps to obtain
$M=H^0_*(C,E\otimes\cO_C(D))$ and to compute
the ideal {\tt J} of $X'\subset\pz^4$,
the embedding of $\pz(E)$ through the linear system $|C_0+p^*(D)|$:
we fix the choice of the curve $C$, named {\tt C}, and 
of the three points $P,Q,Q'$, named respectively 
{\tt p}, {\tt q}, {\tt q'}, and we retype the necessary commands
as in section 4.

The vector space $H^0(\pz^4,\cI_R(3H))$ is therefore given 
by the cubics of $\pz^4$ passing through the fibers over $3Q+3Q'$:
\begin{small}
\begin{verbatim}
betti (qcubeFiber=pullbackIdeal(q^3))
betti (q'cubeFiber=pullbackIdeal(q'^3))
betti (A=intersect(qcubeFiber,q'cubeFiber))
o42 = generators: total: 1 11
                      0: 1  .
                      1: .  .
                      2: . 11
\end{verbatim}
\end{small}

We recall (cf. \S 4) that the ideal {\tt J} of $X'$ contains
a five dimensional space of cubics, and we need to find a set of
representatives cubics {\tt Q} for the quotient space
${H^0(\pz^r,\cI_{R}(3H))}/$ ${H^0(\pz^r,\cI_{X'}(3H))}\cong H^0(X,A)$:
\begin{small}
\begin{verbatim}
betti J
o43 = generators: total: 1 5
                      0: 1 .
                      1: . .
                      2: . 5
a3=super basis(3,A)
j3=super basis(3,J)
Q=super basis (3,ideal a3/ideal j3)
Q=matrix(T,entries Q)
\end{verbatim}
\end{small}

At this point, it remains only to compute the image of $X'$
via the embedding given by $H^0(X,A)$, which is standard.
We report here also the Betti table and the Hilbert function of $X$.
\begin{small}
\begin{verbatim}
Z=K[z_0..z_5];S'=T/J;f=map(S',Z,substitute(Q,S'))
I=ker f --ideal of X
dim I, degree I
o52 = (3, 9)
betti res I
o53 = total: 1 11 18 9 1
          0: 1  .  . . .
          1: .  .  . . .
          2: . 11 18 9 1
apply(0..10,i->hilbertFunction(i,I))
o54 = (1, 6, 21, 45, 78, 120, 171, 231, 300, 378, 465)
\end{verbatim}
\end{small}

\begin{rem}
The projective normality of this surface is proved in the rather
long Proposition 4.6 of \cite{br}.
In \cite{br} they first proved (cf. Lemma 4.3) that 
the projective normality of $X$ is equivalent to 
the 2-normality of $X$, which is equivalent to 
the fact that $X$ does not lies on any quadric.
Thus our computation of the Hilbert function above together
with Lemma 4.3 of \cite{br} verifies the projective normality immediately.
\end{rem}

\end{document}